# ON THE GENERALIZATION OF GURLAND DISTRIBUTION

**Yashwant Singh**

Department of Mathematics,
Seth Moti Lal (P.G.) College, Jhunjhunu(Raj.), India.
E-mail: ysingh23@yahoo.co.in

October 29, 2018

**Abstract**

In the present paper a generalization of Gurland distribution [3] is obtained as a beta mixture of the generalized Poisson distribution (GPD) of Consul and Jain [2]. The first two moments of the distribution and a recurrence relation among probabilities are obtained. The present distribution is supposed to be more general in nature and wider in scope.

**1. Introduction :**

Gurland [3] has obtained a distribution given by its probability mass function (p.m.f.)

$$P(x) = \frac{a(a+1)...(a+x-1)}{(a+b)(a+b+1)....(a+b+x-1)} \phi^x \, _1F_1(a+x; a+b+x, -\phi); x = 0, 1, ...., \quad (1.1)$$

by compounding the Poisson distribution with the beta distribution of first kind. That is, he has considered

$$\text{Poisson } (\theta) \,_{\hat{\theta}/\phi = P} Beta(a, b). \quad (1.2)$$

Here $_1F_1(a;c;x)$ represents the confluent hypergeometric series given by

$$_1F_1(a;c;x) = 1 + \frac{a}{1.c}x + \frac{a(a+1)}{1.2.c(c+1)}x^2 + ... \quad (1.3)$$

The distribution (1.1) was derived by supposing that the number of insect larvae per egg mass has a Poisson distribution with parameter $\theta = \phi p$, where p, which is the probability that an egg hatched into a larva, is assumed to be a random variable having a Beta distribution. The distribution was subsequently studied by Katti [4] who called it type$H_1$ distribution.

The mean and the variance of this distribution are:

$$\mu'_1 = \frac{a\phi}{a+b} \quad (1.4)$$

$$\mu_2 = \frac{a\phi}{a+b} + \frac{ab\phi^2}{(a+b)^2(a+b+1)} \quad (1.5)$$

A generalized verson of the Gurland distribution (1.1) can be obtained using generalized Poisson distribution (GPD) of Consul and Jain [2] given by its pmf

$$P(x) = \frac{\lambda_1(\lambda_1 + x\lambda_2 + x^2\lambda_3)^{x-1} e^{-(\lambda_1 + x\lambda_2 + x^2\lambda_3)}}{x!} \quad (1.6)$$

$\lambda_1 > 0, |\lambda_2| < 1, |\lambda_3| < 1$; x=0,1,2,...

instead of Poisson distribution in (1.2). It can be seen that the Poisson distribution is a particular case of the generalised Poisson distribution just mentioned when $\lambda_2 = 0 = \lambda_3$. The mean and variance of this gneralised Gurland distribution can be obtained as

$$\mu'_1 = \frac{\lambda_1}{(1-\lambda_2-\lambda_3)}, \qquad \mu_2 = \frac{\lambda_1}{(1-\lambda_2-\lambda_3)^3} \qquad (1.7)$$

As the GPD (1.6) is much general in nature and wider in scope, (see Consul [1]) the obtained generalized Gurland distribution is potentially more general in nature and wider in scope.

## 2. A Generalized Gurland Distribution:

The GPD (1.6) can be put in the form

$$P(x) = \frac{\alpha^x (1+x\theta+x^2\phi)^{x-1} e^{-\alpha(1+x\theta+x^2\phi)}}{x!} \qquad (2.1)$$

by putting $\lambda_1 = \alpha$, $\frac{\lambda_2}{\lambda_1} = \theta$, $\frac{\lambda_3}{\lambda_1} = \phi$

We compound this distribution with the beta distribution of first kind in the following way :

$$\text{GPD } (\alpha, \theta, \phi)_{\hat{\theta}/\phi=P} \text{ Beta}(a,b) \qquad (2.2)$$

Thus we find

$$P(x) = \int_0^1 \frac{\alpha^x(1+x\theta+x^2\phi)^{x-1} e^{-\alpha(1+x\theta+x^2\phi)}}{x!} \cdot \frac{1}{B(a,b)} p^{a-1}(1-p)^{b-1} dp \text{ (here } B(a,b) \text{ is the beta function)}$$

$$= \frac{(1+x\theta+x^2\phi)^{x-1}}{x! B(a,b)} \int_0^1 (\delta p)^x e^{-\delta p(1+x\theta+x^2\phi)} p^{a-1}(1-p)^{b-1} dp$$

$$= \frac{\delta^x (1+x\theta+x^2\phi)^{x-1}}{x! B(a,b)} \int_0^1 \Sigma_0^\infty \frac{[-\delta(1+x\theta+x^2\phi)]^s}{s!} p^{a+x+s-1}(1-p)^{b-1} dp$$

$$= \frac{\delta^x (1+x\theta+x^2\phi)^{x-1}}{x! B(a,b)} \Sigma_0^\infty \frac{[-\delta(1+x\theta+x^2\phi)]^s}{s!} B(a+x+s, b)$$

Which after some simplification becomes

$$P(x) = \frac{\delta^x (1+x\theta+x^2\phi)^{x-1}}{x!} \cdot \frac{a(a+1)...(a+x-1)}{(a+b)(a+b+1)...(a+b+x-1)} \cdot {}_1F_1(a+x; a+b+x; -\delta(1+x\theta+x^2\phi)) \qquad (2.3)$$

$x = 0, 1, 2, ...$

The distribution may be termed as the generalized Gurland distribution (GGD).

## 3. Moments :

The mean of GGD (2.3) can be obtained as

$E(X) = E(E(X/P))$

$E(X/P)$ i.e. the conditional expectation of X given P can be obtained by taking $\lambda_1 = \phi p$, $\lambda_2 = \phi\theta p$ and $\lambda_3 = \phi\theta\delta p$ in the expression for mean of the GPD given in (1.7) as

$E(X/P) = \frac{\phi p}{(1-\phi\theta p-\phi\theta\delta p)}$ and thus

$E(X) = E(\frac{\phi p}{1-\phi\theta p-\phi\theta\delta p})$

$= \int_0^1 \frac{\phi p}{1-\phi\theta p-\phi\theta\delta p} \frac{1}{B(a,b)} p^{a-1}(1-p)^{b-1} dp$

$$= \frac{\phi}{B(a,b)} \int_0^1 \Sigma_{s=0}^\infty \frac{(\phi\theta p + \phi\theta\delta p)^s}{s!} p^{a-1}(1-p)^{b-1} dp$$

$$= \phi \Sigma_{s=0}^\infty \frac{B(a+s+1,b)}{B(a,b)}(\phi\theta + \phi\theta\delta)^s$$

After a little simplification thus we find the mean of the GGD as
$$\mu_1' = \frac{\phi a}{(a+b)} \,_2F_1(a+1, a+b+1; \theta\phi + \theta\phi\delta) \qquad (3.1)$$

Where $_2F_1(a,b;c;x)$ represents the Gaussion hypergeometric function given by
$_2F_1(a,b;c;x) = 1 + \frac{a.b}{1.c}x + \frac{a(a+1)b(b+1)}{1.2.c(c+1)}x^2 + ... (3.2)$
Similarly the second moment about origin of the GGD can be obtained as
$E(X^2) = E[E(\frac{X^2}{P})]$
$E(\frac{X^2}{P})$ can be obtained by putting $\lambda_1 = \phi p, \lambda_2 = \phi\theta p$ and $\lambda_3 = \phi\theta\delta p$

in the following expression for $\mu_2'$ obtained from (1.7) :
$$\mu_2' = \frac{\lambda_1}{(1-\lambda_2-\lambda_3)^3} + \frac{\lambda_1^2}{(1-\lambda_2-\lambda_3)^2} \qquad (3.3)$$

as $E(X^2/P) = \frac{\phi p}{(1-\phi\theta p - \phi\theta\delta p)^3} + \frac{\phi^2 p^2}{(1-\phi\theta p - \phi\theta\delta p)^2} \qquad (3.4)$

and thus
$$E(X) = \int_0^1 [\phi p[1 - (\phi\theta p + \phi\theta\delta p)]^{-3} + \phi^2 p^2[1 - (\phi\theta p + \phi\theta\delta p)]^{-2}] \frac{1}{B(a,b)} p^{a-1}(1-p)^{b-1} dp$$

$$= \int_0^1 \phi p[1 + 3(\phi\theta p + \phi\theta\delta p) + 6(\phi\theta p + \phi\theta\delta p)^2 + ...] \frac{1}{B(a,b)} p^{a-1}(1-p)^{b-1} dp$$

$$+ \int_0^1 \phi^2 p^2[1 + 2(\phi\theta p + \phi\theta\delta p) + 3(\phi\theta p + \phi\theta\delta p)^2 + ...] \frac{1}{B(a,b)} p^{a-1}(1-p)^{b-1} dp$$

$$= \frac{\phi}{B(a,b)} \Sigma_{s=1}^\infty \frac{s(s+1)}{2} (\phi\theta + \phi\theta\delta)^{s-1} B(a+s,b) + \frac{\phi^2}{B(a,b)} \Sigma_{s=1}^\infty s(\phi\theta + \phi\theta\delta)^{s-1} B(a+s+1,b)$$

$$= \phi \Sigma_{s=1}^\infty \frac{s(s+1)}{2}(\phi\theta + \phi\theta\delta)^{s-1} \frac{a(a+1)...(a+s+1)}{(a+b)(a+b+1)...(a+b+s)}$$

$$+ \phi^2 \Sigma_{s=1}^\infty s(\phi\theta + \phi\theta\delta)^{s-1} \frac{a(a+1)...(a+s+1)}{(a+b)(a+b+1)...(a+b+s)}$$

$$= \phi \Sigma_{s=1}^\infty s(\phi\theta + \phi\theta\delta)^{s-1} \frac{a(a+1)...(a+s-1)}{(a+b)(a+b+1)...(a+b+s-1)} [\frac{s+1}{2} + \frac{(a+s)}{(a+b+s)}\phi]$$

It can be seen easily that at $\phi = 0 = \theta$, the two moments of the GGD reduce to the respective moments of the Gurland distribution.

**4. Recurrence Relation** :
Denoting the probability function of the GGD (2.1) by $P(x; \phi, \delta, a, b, \theta)$, we have

$$P(x+1;\phi,\delta,a,b,\theta) = \frac{\delta[\delta+(x+1)\delta\theta+(x+1^2\delta\theta]^x}{(x+1)!} \cdot \frac{a(a+1)...(a+x)}{(a+b)(a+b+1)...(a+b+x)}$$
$$_1F_1[a+x+1;a+b+x+1;-\delta(1+(x+1)\theta+(x+1)^2\phi)] \tag{4.1}$$

$$= \frac{\delta[\delta+(x+1)\delta\theta+(x+1^2\delta\phi]^{x-1}[\delta+(x+1)\delta\theta+(x+1^2\delta\phi]}{(x+1)!} \cdot \frac{a(a+1)...(a+x)}{(a+b)(a+b+1)...(a+b+x)}$$
$$\cdot _1F_1[a+x+1;a+b+x+1;-\delta(1+(x+1)\theta+(x+1)^2\phi)]$$

$$= \frac{\delta[\delta+(x+1)\delta\theta+(x+1^2\delta\phi]^{x-1}[-\delta+\delta\theta+\delta\phi]}{(x+1)!} \cdot \frac{a(a+1)...(a+x)}{(a+b)(a+b+1)...(a+b+x)}$$
$$\cdot _1F_1[a+x+1;a+b+x+1;-\delta(1+(x+1)\theta\theta+(x+1)^2\delta\phi)]$$

$$+ \frac{\delta[\delta+(x+1)\delta\theta+(x+1^2\delta\phi]^{x-1}[x\delta\theta+2x\delta\phi]}{(x+1)!} \cdot \frac{a(a+1)...(a+x)}{(a+b)(a+b+1)...(a+b+x)}$$
$$\cdot _1F_1[a+x+1;a+b+x+1;-\delta(1+(x+1)\theta\theta+(x+1)^2\delta\phi)]$$

$$+ \frac{\delta[\delta+(x+1)\delta\theta+(x+1^2\delta\phi]^{x-1}x^2\delta\phi}{(x+1)!} \frac{a(a+1)...(a+x)}{(a+b)(a+b+1)...(a+b+x)}$$
$$\cdot _1F_1[a+x+1;a+b+x+1;-\delta(1+(x+1)\theta\theta+(x+1)^2\delta\phi)]$$

$$= \frac{\delta}{(x+1)}\frac{a}{(a+b)}P(x;\theta+\phi,a+1,b,\theta) + \frac{x}{(x+1)}\frac{\delta(\theta+2\phi)}{(1+\phi+\theta)}P(x;\theta+\phi,a+1,b,\theta)\frac{a}{(a+1)}$$
$$+ \frac{x^2}{(x+1)}\frac{\delta\phi}{(1+\phi+\theta)}P(x;\theta+\phi,a+1,b,\theta)\frac{a}{(a+1)}$$

$$= \frac{\delta}{(x+1)}\frac{a}{(a+b)}P(x;\theta+\phi,a+1,b,\theta)[1 + \frac{x(\theta+2\phi)}{(1+\phi+\theta)} + \frac{x^2\phi}{(1+\phi+\theta)}]$$
$$= \frac{\delta}{(x+1)}\frac{a}{(a+b)}(\frac{(1+\phi+\theta)+x(\theta+2\phi)+x^2\phi}{(1+\phi+\theta)}P(x;\theta+\phi,a+1,b,\theta) \tag{4.2}$$

This recurrence relation among probabilities of the GGD may be helpful in evaluating the probabilities for higher values on the basis of the probabilities for lower values .

**References**:
1. Consul,P.C.;Generalized Poisson distribution , New York , Dekker.
2. Consul,P.C. and Jain,G.C.; A generalization of Poisson distribution, Technometrics,vol.15,791-799.
3.Gurland,J.;A generalized class of contagious distribution, Biometrics,vol.14,229-249.
4.Katti,S.K.;Interrelations among generalized distributions and their components,Biometrics,vol.22,44-52.